\documentclass[english]{article}
\usepackage{amsmath}
\usepackage{amssymb}
\usepackage{amsthm}
\usepackage{babel}

\newtheorem{theo}{Theorem}
\author{Antonio M. Oller Marc\'{e}n}
\title{A Note on Primes Dividing Alternating Sums}
\date{}

\begin{document}
\maketitle

We are all familiar with the harmonic sum:
$$S_n=\sum_{i=1}^n\frac{1}{i}$$
which can easily be shown not to be an integer unless $n=1$ (see
\cite{APO}, chapter 1, exercise 30).

On the other hand we can consider the alternating sum:
$$A_n=\sum_{i=1}^n(-1)^{i-1}\frac{1}{i}=1-\frac12+\frac13-\frac14+\dots+(-1)^{n-1}\frac{1}{n}$$
which is not an integer either. Nevertheless we will be able to
write $\displaystyle{A_n=\frac{a}{b}}$ with $a$ and $b$ being
coprime integers.

Let us start by considering an odd integer $n$ such that
$\displaystyle{\frac{3n+1}{2}=p}$ is a prime. Note that in such a
case it is easily seen that $n$ must satisfy $n\equiv 3\
(\textrm{mod}\ 4)$. Now we construct the sum
$A_n=\displaystyle{\frac{a}{b}}$ and we claim $p$ divides $a$.

In fact we can write
$$A_n=S_n-2\left(\sum_{i=1}^{\frac{n-1}{4}}\frac{1}{2i}\right)=S_n-S_{\frac{n-1}{2}}=\frac{1}{\frac{n-1}{2}+1}+\frac{1}{\frac{n-1}{2}+2}+\dots+\frac{1}{n}$$
Of course, as $p$ is a prime bigger than $n$, we see that the
numbers $\displaystyle{\frac{n-1}{2}+k}$ are, all of them, units in
$\mathbb{Z}/p\mathbb{Z}$ (observe that $k=1,\dots,\frac{n+1}{2}$).
Moreover, there is an even number of summands in
$S_n-S_{\frac{n-1}{2}}$. Now, if we work modulo $p$ and we choose
any $k\in\{1,\dots,\frac{n+1}{2}\}$ we find that
$\displaystyle{\frac{n-1}{2}+k+n-k+1=p}$ so
$\displaystyle{\frac{1}{\frac{n-1}{2}+k}\equiv \frac{-1}{n-k+1}\
(\textrm{mod}\ p)}$ and we get that $A_n\equiv 0\ (\textrm{mod}\ p)$
as claimed.

Now, if we choose an even number $n$ such as
$\displaystyle{\frac{3n+2}{2}=p}$ is a prime, we can reason in the
same way and conclude that $A_n\equiv 0\ (\textrm{mod}\ p)$. Note
that, in this case, it must be $n\equiv 0\ (\textrm{mod}\ 4)$

Finally, we may reformulate the preceding results in the following
way:

\begin{theo}
Let $p$ be and odd prime. Then there exist an integer $n$ such that
$A_n\equiv 0\ (\textrm{mod}\ p)$.
\end{theo}
\begin{proof}
Given an odd prime $p$ it is easy to see that the number $2p-1$ must
be $2p-1\equiv 0,1\ (\textrm{mod}\ 3)$, i.e., it must be $2p-1=3n$
or $2p-1=3n+1$ so it is enough to consider the corresponding $A_n$.
\end{proof}

\end{document}